\documentclass{amsart}
\usepackage{amsfonts}
\usepackage{amscd}
\usepackage{amssymb}
\usepackage{mathrsfs}
\theoremstyle{plain}
\newtheorem{thm}{Theorem}[section]

\theoremstyle{definition}
\newtheorem{ex}[thm]{Example}
\newtheorem{rmk}[thm]{Remark}

\usepackage{amsmath}
\usepackage{amssymb}
\usepackage{amsthm}
\usepackage{latexsym}
\usepackage[T1]{fontenc}
\usepackage[all]{xy}
 \DeclareMathOperator{\fchar}{char}
\DeclareMathOperator{\End}{End} 
\DeclareMathOperator{\Ker}{Ker}

\newcommand{\Cen}{\mathrm{Cen}}

\title[Algebraic dependence of commuting elements
in algebras]{Algebraic dependence of commuting
elements in algebras}

\author{Sergei Silvestrov}
\address{Centre for Mathematical Sciences,
Lund University, Box 118, SE-221 00 Lund, Sweden}
\email{Sergei.Silvestrov@math.lth.se}
\author{Christian Svensson}
\address{Mathematical Institute, Leiden University,
P.O. Box 9512, 2300 RA Leiden, The Netherlands,
and Centre for Mathematical Sciences, Lund
University, Box 118, SE-221 00 Lund, Sweden}
\email{chriss@math.leidenuniv.nl}
\author{Marcel de Jeu}
\address{Mathematical Institute,
Leiden University, P.O. Box 9512, 2300 RA Leiden,
The Netherlands}
\email{mdejeu@math.leidenuniv.nl}



\thanks{This work was supported by a visitor's grant
of the Netherlands Organisation for Scientific
Research (NWO), the Swedish Foundation for
International Cooperation in Research and Higher
Education (STINT), the Crafoord Foundation, the
Royal Physiographic Society in Lund, and the
Royal Swedish Academy of Sciences. We are also
grateful to Lars Hellstr\"om and Daniel Larsson
for helpful comments and discussions.}

\begin{document}

\begin{abstract}
The aim of this paper to draw attention to several aspects of the algebraic
dependence in algebras. The article starts with discussions of the algebraic
dependence problem in commutative algebras. Then the Burchnall-Chaundy
construction for proving algebraic dependence and obtaining the corresponding
algebraic curves for commuting differential operators in the Heisenberg algebra
is reviewed. Next some old and new results on algebraic dependence of commuting
$q$-difference operators and elements in $q$-deformed Heisenberg algebras are
reviewed. The main ideas and essence of two proofs of this are reviewed and
compared. One is the algorithmic dimension growth existence proof. The other is
the recent proof extending the Burchnall-Chaundy approach from differential
operators and the Heisenberg algebra to the $q$-deformed Heisenberg algebra,
showing that the Burchnall-Chaundy eliminant construction indeed provides
annihilating curves for commuting elements in the $q$-deformed Heisenberg
algebras for $q$ not a root of unity.
\end{abstract}

\maketitle

\section{Introduction}
In 1994, one of the authors of the present paper, S.\ Silvestrov, based on
consideration of the previous literature and a series of trial computations,
made the following three part conjecture.

\noindent  % {\bf Conjecture}
\begin{itemize}
\item[$\bullet$] The first part of the conjecture stated that the
Burchnall--Chaundy type result on algebraic dependence of commuting elements
can be proved in greater generality, that is for much more general classes of
non-commutative algebras and rings than the Heisenberg algebra and related
algebras of differential operators treated by Burchnall and Chaundy and in
subsequent literature. \item[$\bullet$] The second part stated that the
Burchnall--Chaundy eliminant construction of annihilating algebraic curves
formulated in determinant (resultant) form works well after some appropriate
modifications for the most or possibly for all classes of algebras where the
Burchnall--Chaundy type result on algebraic dependence of commuting elements
can be proved. \item[$\bullet$] Finally, the third part of the conjecture
stated that the construction and the proof of the vanishing of the
corresponding determinant algebraic curves on the commuting elements can be
performed for all classes of algebras or rings where this fact is true, using
only the internal structure and calculations for the elements in the
corresponding algebras or rings and the algebraic combinatorial expansion
formulas and methods for the corresponding determinants, that is, without any
need of passing to operator representations and use of analytic methods as used
in the Burchnall--Chaundy type proofs.
\end{itemize}

This third part of the conjecture remains widely open with no general such
proofs available for any classes of algebras and rings, even in the case of the
usual Heisenberg algebra and differential operators, and with only a series of
examples calculated for the Heisenberg algebra, $q$-Heisenberg algebra and some
more general algebras, all supporting the conjecture. In the first and the
second part of the conjecture progress has been made. In \cite{HSbook}, the key
Burchnall--Chaundy type theorem on algebraic dependence of commuting elements
in $q$-deformed Heisenberg algebras (and thus as a corollary for $q$-difference
operators as operators representing $q$-deformed Heisenberg algebras) was
obtained. The result and the methods have been extended to more general
algebras and rings generalizing $q$-deformed Heisenberg algebras (generalized
Weyl structures and graded rings) in \cite{LHSSGWSergpotJAlg}. The proof in
\cite{HSbook} is totally different from the Burchnall--Chaundy type proof. It
is an existence argument based only on the intrinsic properties of the elements
and internal structure of $q$-deformed Heisenberg algebras, thus supporting the
first part of the conjecture. It can be used successfully for an algorithmic
implementation for computing the corresponding algebraic curves for given
commuting elements. However, it does not give any specific information on the
structure or properties of such algebraic curves or any general formulae. It is
thus important to have a way of describing such algebraic curves by some
explicit formulae, as for example those obtained using the Burchnall--Chaundy
eliminant construction for the $q=1$ case, i.e., for the classical Heisenberg
algebra. In \cite{LarssonSilv}, a step in that direction was made by offering a
number of examples all supporting the claim that the eliminant determinant
method should work in the general case. However, no general proof for this was
provided. The complete proof following the footsteps of the Burchnall-Chaundy
approach in the case of $q$ not a root of unity has been recently obtained
\cite{JSSvBCqHeis}, by showing that the determinant eliminant construction,
properly adjusted for the $q$-deformed Heisenberg algebras, gives annihilating
curves for commuting elements in the $q$-deformed Heisenberg algebra when $q$
is not a root of unity, thus confirming the second part of the conjecture for
these algebras. In our proof we adapt the Burchnall-Chaundy eliminant
determinant method of the case $q=1$ of differential operators to the
$q$-deformed case, after passing to a specific faithful representation of the
$q$-deformed Heisenberg algebra on Laurent series and then performing a
detailed analysis of the kernels of arbitrary operators in the image of this
representation. While exploring the determinant eliminant construction of the
annihilating curves, we also obtain some further information on such curves and
some other results on dimensions and bases in the eigenspaces of the
$q$-difference operators in the image of the chosen representation of the
$q$-deformed Heisenberg algebra. In the case of $q$ being a root of unity the
algebraic dependence of commuting elements holds only over the center of the
$q$-deformed Heisenberg algebra \cite{HSbook}, and it is unknown yet how to
modify the eliminant determinant construction to yield annihilating curves for
this case.

The present article starts with discussions of the algebraic dependence problem
in algebras. Then the Burchnall-Chaundy construction for proving algebraic
dependence and obtaining corresponding algebraic curves for commuting
differential operators and commuting elements in the Heisenberg algebra is
reviewed. Next some old and new results on algebraic dependence of commuting
$q$-difference operators and elements in the $q$-deformed Heisenberg algebra
are discussed. In the final two subsections we review two proofs for algebraic
dependence of commuting elements in the $q$-deformed Heisenberg algebra. The
first one is the recent proof from \cite{JSSvBCqHeis} extending the
Burchnall-Chaundy approach from differential operators and the Heisenberg
algebra to $q$-difference operators and the $q$-deformed Heisenberg algebra,
showing that the Burchnall-Chaundy eliminant construction indeed provides
annihilating curves for $q$-difference operators and for commuting elements in
$q$-deformed Heisenberg algebras for $q$ not a root of unity. The second one is
the algorithmic dimension growth existence proof from \cite{HSbook}.

\section{description of the problem: commuting
elements in an algebra are given, then find
curves they lie on}

Any two elements $\alpha$ and $\beta$ in a field
$k$ lie on an algebraic curve of the second
degree $F(x,y)=(x-\alpha)(y-\beta) = xy - \alpha
y - \beta x + \alpha \beta = 0.$ The important
feature of this curve is that its coefficients
are also elements in the field $k$. The same
holds if the field $k$ is replaced by any
commutative $k$-algebra $R$, except that then the
coefficients in the {\it annihilating polynomial}
$F$ are elements in $R$, and hence it becomes of
interest from the side of building interplay with
the algebraic geometry over the field $k$ to
determine whether one may find an annihilating
polynomial with coefficients from $k$ for any two
elements in the commutative algebra $R$. It is
well-known that this is not always possible even
in the case of the ordinary commutative
polynomial algebras aver a field unless some
special conditions are imposed on the considered
polynomials. The situation is similar of course
for rational functions or many other commutative
algebras of functions. The ideals of polynomials
annihilating subsets in an algebra are also
well-known to be fundamental for algebraic
geometry, for Gr{\"o}bner basis analysis in
computational algebra and consequently for
various applications in Physics and Engineering.

Another appearance of this type of problems worth
mentioning in our context comes from Galois
theory and number theory in connection to
algebraic and transcendental field extensions. If
$L$ is a field extension $K \subset L$ of a field
$K$, then using Zorn's lemma one can show that
there always exists a subset of $L$ which is
maximal algebraically independent over $K$. Any
such subset is called transcendence basis of the
field $L$ over the subfield $K$ (or of the
extension $K\subset L$). All transcendence bases
have the same cardinality called the
transcendence degree of a field extension. If
$B=\{b_1,\dots, b_n\}$ is a finite transcendence
basis of $K\subset L$, then $L$ is an algebraic
extension of the subfield $K(B)$ in $L$ generated
by $B$ over the subfield $K$. This means in
particular that for any $l\in L$ the set
$\{B,l\}$ is algebraically dependent over $K$,
that is, there exists a polynomial $F$ in $n+1$
indeterminates over $K$ such that $F(b_1,\dots,
b_n,a) = 0$. Investigation of algebraic
dependence, transcendence basis and transcendence
degree is a highly non-trivial important
direction in number theory and theory of fields
with striking results and many longstanding open
problems. For example, if algebraic numbers
$a_1,...,a_n$ are linearly independent over
$\mathbb{Q}$, then $e^{a_1},...,e^{a_n}$ are
algebraically independent over $\mathbb{Q}$
(Lindemann-Weierstrass theorem, 1880's). Whether
$e$ and $\pi$ are algebraically dependent over
$\mathbb{Q}$ or not is unknown, and only
relatively recently (1996) a long-standing
conjecture on algebraic independence of $\pi$ and
$e^\pi$ was confirmed \cite{Nesterenko}.
Comprehensive overviews and references in this
direction can be found in \cite{NestPhil}.

In view of the above, naturally important
problems are to describe, analyze and classify:
\begin{itemize}
\item[(C1)] commutative algebras over $K$ in which any pair of elements is
algebraically dependent over $K$;
 \item[(C2)] pairs $(A, B)$ of a
commutative algebra $A$ and a subalgebra $B$ over a field $K$ such that any
pair of elements in $A$ is algebraically dependent over $B$.
\end{itemize}
Problem (C1) is of course a special case of (C2).

In the polynomial algebra $K[x_1,\dots,x_n]$
generated by $n$ independent commuting
indeterminates, for instance the set
$\{x_1,\dots,x_n\}$ as well as any of its
non-empty subsets are algebraically independent
over $K$. Thus $K[x_1,\dots,x_n]$ does not belong
to the class of algebras in the problem
\mbox{(C1)}. The same holds of course for any
algebra containing $K[x_1,\dots,x_n]$. In general
algebraic dependence of polynomials happens only
under some restrictive conditions on the rank (or
vanishing) of their Jacobian and other fine
aspects.

This indicates that some principal changes have
to be introduced in order to be able to get
examples of algebras satisfying \mbox{(C1)} or
\mbox{(C2)}. It turns out that the range of
possibilities expands dramatically within the
realm of non-commutative geometry, if the
indeterminates (the coordinates) are
non-commuting.

In any algebra there are always commutative
subalgebras, for instance any subalgebra
generated by any single element. If the algebra
contains for instance a non-zero nilpotent
element ($a^n=0$ for some $n>2$), then the
subalgebra generated by this element satisfies
$(C1)$. Indeed, any two elements $f_1= ap(a)+p_0$
and $f_2= aq(a)+q_0$ in this commutative
subalgebra are annihilated by $F(s,t)=(q_0 s -
p_0 t)^n$, since $F(f_1,f_2) =(q_0 f_1 - p_0
f_2)^n =a^n (q_0 p(a) - p_0 q(a))^n = 0$.

For the non-commutative algebras the problems
inspired by $C1$ and $C2$ are to describe,
analyze and classify:
\begin{itemize}
\item{(NC1)} algebras over a field $K$ such that in all their commutative
subalgebras any pair of elements is algebraically dependent over the field;
\item{(NC2)} pairs $(A, B)$ of an algebra $A$ and a subalgebra $B$ over a field
$K$ such that any pair of commuting elements in $A$ is algebraically dependent
over $B$.
\end{itemize}
Problem (NC1) is of course a special case of
(NC2).

From the situation in commutative algebras as described above, one might
intuitively expect that finding non-commutative algebras satisfying (NC1) is
difficult if not impossible task. However, this perception is not quite
correct. One of general purposes for the present article is to illuminate this
phenomena.

\section{Burchnall-Chaundy construction for differential operators}

Among the longest known, constantly studied and
used all over Mathematics non-commutative
algebras is the so called Heisenberg algebra,
also called Weyl algebra or Heisenberg-Weyl
algebra. It is defined as an algebra over a field
$K$ with two generators $A$ and $B$ and defining
relations $AB-BA=I$, or equivalently as
${\mathcal
H_1}=K\left<A,B\right>/\left<AB-BA-I\right> $,
the quotient algebra of a free algebra on two
generators $A$ and $B$ by the two-sided ideal
generated by $AB-BA-I$ where $I$ is the unit
element.

It's ubiquity is due to the fact that the basic
operators of differential calculus $A= D=
\frac{d}{dx}$ and $M_x: f(x) \mapsto x f(x)$
satisfy the Heisenberg algebra defining relation
$DM_x-M_xD=I$ on polynomials, formal power
series, differentiable functions or any linear
spaces of functions invariant under these
operators due to the Leibniz rule. On any such
space, the representation $(D,M_x)$ of the
commutation relation $AB-BA=I$ defines a
representation of the algebra ${\mathcal H_1}$,
called sometimes, especially in Physics, the
(Heisenberg) canonical representation. We will
use sometimes this term as well, for shortage of
presentation.

If $K$ is a field of characteristic zero, then the algebra $${\mathcal
H_1}=K\left<A,B\right>/\left<AB-BA-I\right> $$ is simple, meaning that it does
not contain any two-sided ideals different from zero ideal and the whole
algebra. The kernel $\Ker (\pi) = \{a \in A \mid \pi (a)=0\}$ of any
representation $\pi$ of an algebra is a two-sided ideal in the algebra. Thus
any non-zero representation of ${\mathcal H_1}$ is faithful, which is in
particular holds also for a canonical representation. Thus the algebra
${\mathcal H_1}$ is actually isomorphic to the ring of differential operators
with polynomial coefficients acting for instance on the linear space of all
polynomials  or on the space of formal power series in a single variable.

In the literature on algebraic dependence of
commuting elements in the Heisenberg algebra and
its generalizations -- a result which is
fundamental for the algebro-geometric method of
constructing and solving certain important
non-linear partial differential equations -- one
can find several different proofs of this fact,
each with its own advantages and disadvantages.
The first proof of such result utilizes
analytical and operator theoretical methods. It
was first discovered  by Burchnall and Chaundy
\cite{Burchnall&Chaundy1} in the 1920's. Their
articles
\cite{Burchnall&Chaundy1,Burchnall&Chaundy2,
Burchnall&Chaundy3} contain also pioneering
results in the direction of in-depth connections
to algebraic geometry. These fundamental papers
were largely forgotten for almost fifty years
when the main results and the method of the
proofs of Burchnall and Chaundy were rediscovered
in the context of integrable systems and
non-linear differential equations
\cite{KrichIntro1,KrichIntro2,MumIntro}. Since
the 1970's, deep connections between algebraic
geometry and solutions of non-linear differential
equations have been revealed, indicating an
enormous richness largely yet to be explored, in
the intersection where non-linear differential
equations and algebraic geometry meet. This
connection is of interest both for its own
theoretical beauty and because non-linear
differential equations appear naturally in a
large variety of applications, thus providing
further external motivation and a source of
inspiration for further research into commutative
subalgebras and their interplay with algebraic
geometry.

A second, more algebraic approach to proving the algebraic dependence of
commuting differential operators was obtained in a different context by Amitsur
\cite{Amitsur} in the 1950's. Amitsur's approach is more in the direction of
the classical connections with field extensions we have already mentioned.
Recently in the 1990's a more algorithmic combinatorial method of proof based
on dimension growth considerations has been found in \cite{HDiscMath,HSbook}.
The main motivating problem for these developments was to describe, as detailed
as possible, algebras of commuting differential operators and their properties.
The solution of this problem is where the interplay with algebraic geometry
enters the scene. The Burchnall-Chaundy result is responsible for this
connection as it states that commuting differential operators satisfy equations
for certain algebraic curves, which can be explicitly calculated for each pair
of commuting operators by the so called eliminant method. The formulas for
these curves, obtained from this method by using the corresponding
determinants, are important for their further analysis and for applications and
further development of the general method and interplay with algebraic
geometry.

In the rest of this section we will briefly
review the basic steps of the Burchnall--Chaundy
construction. For simplicity of exposition until
the end of this section we assume  that the field
of scalars is $K=\mathbb{C}$.

Commutativity of a pair of differential operators
$$P=\sum_{i=0}^{m} p_i(t)\partial^i, \qquad\quad
Q=\sum_{i=0}^{n} q_i(t)\partial^i$$ where
$p_i,q_i$ are analytic functions in $t$ and
$\partial:=\frac{\mathrm{d}}{\mathrm{d}t}$, puts
severe restrictive conditions on the functions
$p_i$ and $q_i$. In its original formulation the
result of Burchnall and Chaundy can be stated as
follows.
\begin{thm}[Burchnall--Chaundy, 1922]
For any two commuting
  differential operators $P$ and $Q$, there is a polynomial
  $F(x,y)\in  \mathbb{C}[x,y]$ such that
  $F(P,Q)=0$.
\end{thm}

The polynomial appearing in this theorem is often
referred to as the \emph{Burchnall--Chaundy
polynomial}.

It is worth mentioning that in their papers Burchnall and Chaundy have neither
specified any conditions on what kind of functions the coefficients in the
differential operators are, nor the spaces on which these operators act. Thus
more precise formulations of the result should contain such a specification.
Any space of functions where the construction is valid would be fine (e.g.,
polynomials or analytic functions in the complex domain). Algebraic steps in
the construction are generic. However in order to reach the main conclusions on
the existence and annihilating property of the curves, the Burchnall-Chaundy
considerations use existence of solutions of an eigenvalue problem for ordinary
differential operators and the property that the dimension of the solution
space of a homogeneous differential equation does not exceed the order of the
operator. To ensure these properties, coefficients of differential operators
must be required to belong to not very restrictive but nevertheless specific
classes of functions.

We will now sketch Burchnall-Chaundy construction for the convenience of the
reader and in connection with further considerations in the present article. In
spite of the fact that the Burchnall-Chaundy arguments actually do not
constitute a complete proof due to some serious gaps, they provide important
insight for building annihilating curves which can be developed into a complete
proof and a well functioning construction after appropriate adjustments and
restructuring.

If differential operators $P$ and $Q$ of orders $m$ and $n$, respectively,
commute then $P-h$ and $Q$ commute for any constant $h\in \mathbb{C}$. Thus
$Q(\Ker (P-h)) \subseteq \Ker (P-h)$. Consequently, if $y_1, \dots, y_m$ is a
basis of $\Ker (P-h)$, the fundamental set of solutions of the eigenvalue
problem for the differential equation $P(y)-hy = 0$ (note that the existence
and the dimension properties of the solution space are assumed to hold here),
then $Q(y_1),\dots,Q(y_m)$ are also elements of $\Ker (P-h)$ and hence
$Q(\vec{y})= A \vec{y}$ for some matrix $A = (a_{i,j})_{i,j=1}^{m}$ with
entries from $\mathbb{C}$. Let $k\in \mathbb{C}$ be another arbitrary constant.
A common nonzero solution, $Y= \vec{c}^T\vec{y}, \vec{c}\in \mathbb{C}^n $, of
eigenvalue problems $PY=hY,\ QY=hY$, or equivalently a nonzero $Y \in \Ker
(P-h) \cap \Ker (Q-k)$, exists if and only if $(Q-k)Y=\vec{c}^{\ T} (A-k)
\vec{y} =0$ has a nonzero solution $\vec{c}$, which happens only if $\det (A-k)
= 0$. This is a polynomial in $k$ of order $m$ and hence corresponding to each
$h$ there exists only $m$ values of the constant $k$ (not necessarily all
distinct) such that there exists nonzero $Y \in \Ker (P-h) \cap \Ker (Q-k)$
Note that here it was used that the scalar field is algebraically closed.
Similarly, corresponding to each $k$ there exists $n$ values of $h$ such that
$\Ker (P-h) \cap \Ker (Q-k) \neq \{0\}$. From this Burchnall and Chaundy
conclude that if $\Ker (P-h) \cap \Ker (Q-k) \neq \{0\}$, then $h$ and $k$
satisfy some polynomial equation $F(h,k)=0$, where $F$ is a polynomial of
degree $n$ in $h$ and $m$ in $k$ with coefficients from $\mathbb{C}$. There is
however a problem with this key conclusion. Surely, the equation $h=2k + \sin
(k)$ gives a bijection between the $h$'s and the $k$'s, but the resulting curve
is not algebraic. Thus as we already pointed out, substantial adjustments are
necessary in order to rearrange these arguments into a functioning construction
and a complete proof. Suppose, however, that the proper adjustments have been
made, and that a polynomial $F$ with the above annihilating properties has been
found. Then, any $Y\in \Ker (P-h) \cap \Ker (Q-k)$ is also a solution of the
differential equation $F(P,Q) Y = F(h,k) Y =0$ which is of order $mn$ unless it
happens that $F(P,Q) = 0$. Thus there can be at most $mn$ linearly independent
$Y\in \Ker (P-h) \cap \Ker (Q-k)$. Note that here again a specific property of
the dimension of the solution space of a differential equation is assumed to
hold. For each $h$ there exists $k$ such that $\Ker (P-h) \cap \Ker (Q-k) \neq
\{0\}$. Since the field $\mathbb C$ is infinite (note that it is another
special property of the field), one can chose infinitely many pairwise distinct
numbers $h$ with corresponding $k$ and nonzero functions $Y_{h,k} \in \Ker
(P-h) \cap \Ker (Q-k)$. But any nonempty set of eigenfunctions with pairwise
distinct eigenvalues for a linear operator is always linearly independent. Thus
the dimension of the solution space $\Ker(F(P,Q))$ is infinite. But this
contradicts to the already proved $\Ker(F(P,Q))\leq mn $ unless $F(P,Q)=0$.
Therefore, indeed $F(P,Q)=0$ which is exactly what was claimed.

A beautiful feature of the Burchnall-Chaundy
arguments in the differential operator case,
however, is that they are almost constructive in
the sense that they actually tell us, after
taking a closer look, how to compute such
annihilating curves, given the commuting
operators. This is done by constructing the
\emph{resultant} (or \emph{eliminant}) of
operators $P$ and $Q$. We sketch this
construction, as it is important to have in mind
for this article. To this end, for
 complex variables $h$ and $k$, one writes:
\begin{align}\label{diff_op_res_rowP}
    \partial^r(P-h\mathbf{1})=\sum_{i=0}^{m+r}
    \theta_{i,r}\partial^i-h\partial^r, \quad\quad
    r=0,1,\dots ,n-1\\
    \label{diff_op_res_rowQ}
    \partial^r(Q-k\mathbf{1})=\sum_{i=0}^{n+r}
    \omega_{i,r}\partial^i-k\partial^r, \quad\quad
    r=0,1,\dots ,m-1
\end{align}
where $\theta_{i,r}$ and $\omega_{i,r}$ are certain functions built from the
coefficients of $P$ and $Q$ respectively, whose exact form is calculated by
moving $\partial^r$ through to the right of the coefficients, using the Leibniz
rule. The coefficients of the powers of $\partial$ on the right hand side in
(\ref{diff_op_res_rowP}) and (\ref{diff_op_res_rowQ}) build up the rows of a
matrix exactly as written. That is, as the first row we take the coefficients
in $\sum_{i=0}^{m} \theta_{i,0}\partial^i-h\partial^0$, and as the second row
the coefficients in $\sum_{i=0}^{m+1} \theta_{i,1}\partial^i-h\partial$,
continuing this until $k=n-1$. As the $n^{\mathrm{th}}$ row we take the
coefficients in $\sum_{i=0}^{n} \omega_{i,0}\partial^i-k\partial^0$, and as the
$(n+1)^{\mathrm{th}}$ row we take the coefficients in $\sum_{i=0}^{n+1}
\omega_{i,1}\partial^i-k\partial$ and so on. In this manner we get a
$(n+m)\times (n+m)$-matrix using (\ref{diff_op_res_rowP}) and
(\ref{diff_op_res_rowQ}). The determinant of this matrix yields a trivariate
polynomial $F(x,h,k)$ over $\mathbb{C}$. When written as $F(x,h,k) = \sum_i
\delta_i (h,k) x^i$, it can be proved, using existence and uniqueness results
for ordinary differential equations, that $\delta_i(P,Q) = 0$ for all $i$. It
is not difficult to see that the $\delta_i$ are not all zero.

For clarity, we include the following example.
\begin{ex}\label{eli32}
Let $P$ and $Q$ be as above, with $m =3$ and $n=2$. We then have \\ $F(x,h,k)=$
\[\left|\begin{array}{ccccc}
p_{0,0}(x) - h & p_{0,1}(x) & p_{0,2} (x) &
p_{0,3} (x) & 0\\
p_{1,0}(x)  & p_{1,1} (x)  - h & p_{1,2} (x)  &
p_{1,3} (x) & p_{1,4} (x)\\
q_{0,0}(x) - k & q_{0,1}(x) & q_{0,2} (x) & 0 & 0\\
q_{1,0}(x)  & q_{1,1} (x)  - k & q_{1,2} (x)  &
q_{1,3} (x) & 0\\
q_{2,0}(x)  & q_{2,1} (x)  & q_{2,2}(x) - k & q_{2,3}(x) & q_{2,4} (x)
\end{array}\right|.\]
\end{ex}

\begin{rmk} Thus the pairs of eigenvalues
$(h,k)\in {\mathbb C}^2$ corresponding to the same joint eigenfunction, i.e.,
corresponding to the same non-zero $y$ such that
\begin{align*}
    Py=hy\qquad\text{and}\qquad Qy&=ky,
\end{align*} lie on curves $Z_i$ defined by the $\delta_i$. There is a ``dictionary'' between certain
geometric and analytic data \cite{MumIntro} allowing in particular construction
of common eigenfunctions for commuting operators using geometric data
associated to the compactification of the annihilating algebraic curve.
\end{rmk}

Burchnall--Chaundy result for operators with
polynomial coefficients can be reformulated in
more general terms for the abstract Heisenberg
algebra $\mathcal{H}_1$ rather than in terms of
its specific canonical representation by
differential operators. This specialization of
coefficients to polynomials does not influence
the Burchnall--Chaundy construction of the
annihilating algebraic curves. Thus restricting
to this context does not effect the main
ingredients needed for building the interplay
with algebraic geometry. At the same time, when
moreover reformulated entirely in the general
terms of the algebra $\mathcal{H}_1$, this
specialization becomes a generalization in
another way, because establishing algebraic
dependence of commuting elements directly in the
Heisenberg algebra, without passing to the
canonical representation, makes the result valid
not just for differential operators, but also for
any other representations of $\mathcal{H}_1$ by
other kinds of operators.

\begin{thm}["Burchnall--Chaundy", algebraic version] Let $P$ and $Q$ be two commuting
  elements in $\mathcal{H}_1$, the Heisenberg algebra. Then there is a
  bivariate polynomial
  $F(x,y)\in {\mathbb C}[x,y]$ such that $F(P,Q)=0$.
\end{thm}

\section{Burchnall--Chaundy theory for the $q$-deformed Heisenberg algebra}

Let $K$ be a field. If $q \in K$ then $\mathcal{H}_q=\mathcal{H}_K (q)$, the
$q$-deformed Heisenberg algebra over $K$, is the unital associative $K$-algebra
which is generated by two elements $A$ and $B$, subject to the commutation
relation $AB - qBA = I$. Though this algebra sometimes is also called the
$q$-deformed Weyl or $q$-deformed Heisenberg-Weyl algebra in various contexts,
we will follow systematically the terminology in \cite{HSbook}. We will
indicate how one can prove that - under a condition on $q$ - for any commuting
$P, Q \in \mathcal{H}_q$ of order at least one (where ``order'' will be defined
below), there exist finitely many explicitly calculable polynomials $p_i \in
K[X, Y]$ such that $p_i (P, Q) =0$ for all $i$, and at least one of the $p_i$
is non-zero. The number of polynomials depends on the coefficients of $P$ and
$Q$, as well as on their orders. The polynomials $p_i$ can be obtained by the
mentioned before so-called eliminant construction which was introduced for
differential operators (the case of $q=1$) by Burchnall and Chaundy in
\cite{Burchnall&Chaundy1,Burchnall&Chaundy2,Burchnall&Chaundy3} in 1920's, and
which we will employ and extend to the context of general $q$-deformed
Heisenberg algebras showing that analogous determinant (resultant) construction
of the annihilating algebraic curves works for $q$-deformed Heisenberg algebras
well.

We assume $q \neq 0$ throughout (our method
breaks down when $q = 0$) and define the
$q$-integer $\{n\}_q$, for $n \in \mathbb{Z}$, by
\[\{n\}_q = \left\{ \begin{array}{lr}
\frac{q^n -1}{q-1} & q \neq 1;\\
 n & q = 1.
\end{array}\right.\]
Following \cite[Definition 5.2]{HSbook} we will
say that $q\in K^*=K\setminus \{0, 1\}$ is of
{\it torsion type} if there is a positive integer
solution $p$ to $q^p=1$, and of {\it free type}
if the only integer solution to $q^p=1$ is $p=0$.
In the torsion type case the least such positive
integer $p$ is called the order of $q$ and in
free type case the order of $q$ is said to be
zero. If $q=1$ then it is said to be of a torsion
type if $K$ is a field of non-zero
characteristic, and of free type if the field is
of characteristic zero.
\\

\begin{rmk}\label{nonunroot}
The following are equivalent for $q \neq 0$:
\nopagebreak
\begin{enumerate}
\item for $n \in \mathbb{Z}$, $\{n\}_q =0$
if and only if $n=0$;
\item for $n_1, n_2 \in \mathbb{Z}$,
$\{n_1\}_{q} = \{n_2\}_{q}$ if and only if $n_1 =
n_2$;
\item $\left\{\begin{array}{lr}
q \textup{ is not a root of unity other than } 1, &
\textup{if } \fchar k = 0;\\
q \textup{ is not a root of unity},  &
\textup{if } \fchar k \neq 0.\\
\end{array}
\right.$
\end{enumerate}

Part (2) of this remark is essential when one
considers the dimension of eigenspaces for
$q$-difference operators which is important for
our extension of the Burchnall-Chaundy
construction.

Note also that under our assumptions $K$ is
infinite.
\end{rmk}

Let $\mathcal{L}$ be the $K$-vector space of all
formal Laurent series in a single variable $t$
with coefficients in $K$. Define
\[M (\sum_{n= -\infty}^{\infty} a_n t^n) = \sum_{n = - \infty}^{\infty} a_n t^{n+1} = \sum_{n = - \infty}^{\infty} a_{n-1} t^n,\]
\[D_q (\sum_{n= -\infty}^{\infty} a_n t^n) = \sum_{n = - \infty}^{\infty} a_n \{n\}_q t^{n-1} = \sum_{n = - \infty}^{\infty} a_{n+1}
\{n+1\}_{q} t^{n}.\] Alternatively, one could
introduce $\mathcal{L}$ as the vector space of
all functions from $\mathbb{Z}$ to $K$ and let
$M$ act as the right shift and $D_q$ as a
weighted left shift, but the Laurent series model
is more appealing.

The algebra $\mathcal{H}_q$ has $\{I, A, A^2,
\ldots\}$ as a free basis in its natural
structure as a left $K[X]$-module with $X$ mapped
to $B$.

So an arbitrary non-zero element $P$ of
$\mathcal{H}_q$ can be written as
\[P = \sum_{j=0}^m p_j (B) A^j, p_m \neq 0,\]
for uniquely determined $p_j \in K[X]$ and $m
\geq 0$. The integer $m$ is called the
\emph{order} of $P$ (or degree with respect  to
A) \cite{HSbook}.

By sending $A$ to $D_q$ and $B$ to $M$,
$\mathcal{L}$ becomes an $\mathcal{H}_q$-module.
If we make the additional assumption that
$\{n\}_q \neq 0$ or equivalently $q^n\neq 1$ for
all non-zero $n \in \mathbb{Z}$, then this
representation is faithful \cite[Theorem
8.3]{HSbook}.

We will assume that $q \neq 0$ and $\{n\}_q \neq
0$ for $n \neq 0$ throughout this paper and
identify $\mathcal{H}_q$ with its image in
$\End_K(\mathcal{L})$ under the previously
defined representation without further notice. In
the image of $\mathcal{H}_q$ under the
representation, $\{1, D_q, D_q^2, \ldots\}$ is a
free basis in its natural structure as a left
$K[X]$-module with $X$ being mapped to $M$, and
so any $P$ can be written as
\[P = \sum_{j=0}^m p_j (M) D_q^j,
\ p_m \neq 0,\] for uniquely determined $p_j \in
K[X]$ and $m \geq 0$. The integer $m$ is called
the \emph{order} of $P$.

The important result in the context of this
article is an extension of the Burchnall--Chaundy
theorem in algebraic version for the $q$-deformed
Heisenberg algebra $\mathcal{H}_q$, due to
Hellstr\"om and Silvestrov \cite[Therem
7.5]{HSbook}.
\begin{thm}[Hellstr\"om--Silvestrov,
\cite{HSbook}] \label{thm:qBCHS} Let $q\in
K^*=K\setminus \{0\}$, and let $P$ and $Q$ be two
commuting
  elements in $\mathcal{H}_q$. Then there is a
  bivariate polynomial
  $F(x,y)\in\mathrm{Z}(\mathcal{H}_q)[x,y]$,
  with coefficients in the center
  $\mathrm{Z}(\mathcal{H}_q)$
  of $\,\mathcal{H}_q$, such that $F(P,Q)=0$.
\end{thm}
If $q$ is not a root of unity ($\{n\}_q$ are different for different $n$), the
center of $\mathcal{H}_q$ is trivial, i.e., consists of scalar multiples of the
identity $\mathrm{Z}(\mathcal{H}_q)=K I$. Thus in this case there exists a
"genuine" algebraic curve over the scalar field $K$ as Theorem \ref{thm:qBCHS}
takes the following form \cite[Theorem 7.4]{HSbook}.
\begin{thm}[Hellstr\"om--Silvestrov,
\cite{HSbook}] \label{thm:qBCHSqfreetype} Let
$q\in K^*=K\setminus \{0\}$ be of free type, and
let $P$ and $Q$ be two commuting elements in
$\mathcal{H}_q$. Then there is a bivariate
polynomial $F(x,y)\in K[x,y]$, with coefficients
in $K$, such that $$F(P,Q)=0.$$
\end{thm}
Note however, that when $q$ is of torsion type and order $d$ (i.e., $d$ is the
smallest positive integer such that $q^d=1$), then the center is
$$\mathrm{Z}(\mathcal{H}_q)=K[A^d,B^d],$$
the subalgebra spanned by $\{A^d,B^d\}$, where
$d$ is the order of $q$, the minimal positive
integer such that $q^d=1$ (\cite[Corollary
6.12]{HSbook}). The conclusion of
Theorem~\ref{thm:qBCHSqfreetype}, that is the
algebraic dependence of commuting elements over
the field $K$, does not hold for $q$ of torsion
type, since if $p$ is the order of $q$, then
\(\alpha=A^p\) and \(\beta=B^p\) commute, but do
not satisfy any commutative polynomial relation.
Thus in this case Theorem \ref{thm:qBCHS} has
indeed to be invoked.

The proof of Theorem \ref{thm:qBCHS} as given in \cite{HSbook} is purely
existential. However, while it says essentially nothing theoretically on the
form or properties of the annihilating curves, the construction used in the
proof actually provides an explicit computationally implementable algorithm for
producing annihilating  polynomials.

A specialization of a part of the general conjecture made by S.\ Silvestrov in
1994 is that the determinant scheme devised by Burchnall and Chaundy could be
used to calculate the polynomial even in the case of $\mathcal{H}_q$. The
Burchnall-Chaundy eliminant construction adaptation to $\mathcal{H}_q$ and a
series of examples, indicating that the construction indeed yields annihilating
algebraic curves for commuting elements in $\mathcal{H}_q$, was first presented
in \cite{LarssonSilv}. But no full proof of the conjecture that this adaption
is possible for all commuting elements of $\mathcal{H}_q$ was given. Even
though a direct generalization of the classical arguments of Burchnall and
Chaundy was attempted some important technical steps were missing. The main
reason why an analogous proof for $q$-difference operators (i.e., elements in
$\mathcal{H}_q$) is problematic is that the solution space for $q$-difference
equations may not in general be as well-behaved as for ordinary differential
operators, and the proof of the classical Burchnall--Chaundy theorem relies
heavily on properties of the solution spaces to the eigenvalue-problems for the
differential operators $P$ and $Q$.

Nevertheless, in the case of $q$ of free type we have succeeded now to extend
the original Burchnall-Chaundy eliminant method to the $q$-difference operators
and hence to $\mathcal{H}_q$. Roughly speaking, the key technical idea is to
choose an appropriate representation space for the canonical representation of
$\mathcal{H}_q$ generated by $M_x, D_q$,  so that the necessary key ingredients
about dimensions of eigenspaces of $q$-difference operators in the
Burchnall-Chaundy type proof are still available. The module $\mathcal{L}$
introduced above has all these properties. It is thus that, for $q$ of free
type, we again have a determinant based construction providing annihilating
curves. The complete proofs will be presented in \cite{JSSvBCqHeis}; they are
considerably more involved than the original work by Burchnall and Chaundy.

An extension of the result and the Burchnall-Chaundy type construction to the
case of $q$ of torsion type, i.e., when $\mathcal{H}_q$ has non-trivial center,
is still not available. This is just one of many reasons why another proof of
the possibility of adapting to $\mathcal{H}_q$ the determinant construction of
the annihilating curves, relying on purely algebraic methods, i.e., without
passing to a specific representation of $\mathcal{H}_q$, is desirable. The
existence of such proof is a specialization to $\mathcal{H}_q$ of the third
part of the aforementioned general conjecture by S.\ Silvestrov.

\subsection{Eliminant determinant construction
for $q$-deformed Heisenberg algebra}

In this section we will briefly outline our
extension of the Burchnall-Chaundy result and
constructions to the $q$-difference operators, or
to the abstract algebra context of the
$q$-deformed Heisenberg algebra. The complete
details of the proofs will be presented in
\cite{JSSvBCqHeis}.

Let $P, Q \in \mathcal{H}_q$ be of order $m \geq
1$ and $n \geq 1$ respectively. Write for $k = 0,
\ldots, n-1$,
\[D_q^k P = \sum_{j=0}^{m+k} p_{k, j} (M) D_q^j, \textup{ with } p_{k,j} \in K[X],\]
and, for $l = 0, \ldots, m-1$, write
\[D_q^l Q = \sum_{j=0}^{n+l} q_{l, j} (M) D_q^j,
\textup{ with } q_{l,j} \in K[X].\] By analogy with the Burchnall-Chaundy
method for differential operators, we build up an $(m+n) \times (m+n)$-matrix
as follows. For $k = 1, \ldots, n$, the $k$-th row is given by the coefficients
of powers of $D_q$ in the expression $D_q^{k-1} P - \lambda D_q^{k-1}=
\sum_{j=0}^{m+k-1} p_{k-1, j} (M) D_q^j - \lambda D_q^{k-1}$, where $\lambda$
is a formal variable. For $k \in \{ n+1, \ldots, m+n\}$, the $k$-th row is
given by the coefficients of $D_q$ in $D_q^{k-n-1} Q - \mu D_q^{k-n-1}=
\sum_{j=0}^{k-1} p_{k-n-1, j} (M) D_q^j - \mu D_q^{k-n-1}$, where $\mu$ is a
formal variable different from $\lambda$. The determinant of this matrix is
called the \emph{eliminant} of $P$ and $Q$. We denote it $\Delta_{(P, Q)} (M,
\lambda, \mu)$. It is a polynomial with coefficients in $K$.

\begin{thm} {\rm (\cite{JSSvBCqHeis})}
\label{mainthm} Let $K$ be a field and $0 \neq q
\in K$ be such that $\{n\}_q =0$ if and only if
$n=0$. Suppose
\[P = \sum_{j=0}^m p_j (M) D_q^j \,\,\, (m \geq 1, p_m \neq 0)\]
and
\[Q = \sum_{j=0}^n q_j (M) D_q^j \,\,\, (n \geq 1, p_n \neq 0)\]
are commuting elements of $\mathcal{H}_q$, and
let $\Delta_{P,Q} (M, \lambda, \mu)$ be the
eliminant constructed as above. Then
$\Delta_{P,Q} \neq 0$. In fact, if $\,\,q_n (M) =
\sum_{i} a_i M^i \,\,(a_i \in K)$ then
$\Delta_{P,Q}$ has degree $n$ as a polynomial in
$\lambda$, and its non-zero coefficient of
$\lambda^n$ is equal to $(-1)^n \prod_{k=0}^{m-1}
(\sum_i a_i q^{k\cdot i} M^i)$. Likewise, if $p_m
(M) = \sum_i b_i M^i$, then $\Delta_{P,Q}$ has
degree $m$ as a polynomial in $\mu$, and its
non-zero coefficient of $\mu^m$ is equal to
$(-1)^m \prod_{k=0}^{n-1} (\sum_i b_i q^{k \cdot
i} M^i)$.

Let \begin{displaymath}s = n \cdot \max_{j}
\deg(p_j) +  m \cdot\max_j
\deg(q_j),\end{displaymath}
\[t = \frac{1}{2} \cdot n\cdot (n-1) \cdot \max_j \deg(p_j) + \frac{1}{2} \cdot m\cdot  (m-1)\cdot \max_j \deg(q_j),\]
and define the polynomials $\delta_i$ $(i = 1,
\ldots, s)$ by $\Delta_{P,Q} (M, \lambda, \mu) =
\sum_{i=0}^s \delta_i (\lambda, \mu) M^i$.

Then

\begin{enumerate}
\item each of the coefficients of the $\delta_i$
can be expressed as $\sum_{l=0}^t r_l q^l$ with
the $r_l$ in the subring of $K$ which is generated
by the coefficients of all the $p_j$ and the $q_j$,
\item at least one of the $\delta_i$ is non-zero,
\item $\delta_i (P, Q) = 0$
for all $i = 0, \ldots, s$.
\end{enumerate}
\end{thm}
The reader will easily convince himself of all
statements in the theorem other than (3). We will
now sketch the main idea of the proof of (3) as
given in \cite{JSSvBCqHeis}.

The idea is as follows. Suppose $\lambda_0, \mu_0
\in K$ and $0 \neq v_{(\lambda_0, \mu_0)} \in
\mathcal{L}$ is a common eigenvector of $P$ and
$Q$:
\[P v_{(\lambda_0, \mu_0)} = \lambda_0
v_{(\lambda_0, \mu_0)},\]
\[Q v_{(\lambda_0, \mu_0)} = \mu_0
v_{(\lambda_0, \mu_0)}.\] Then the specialization
$\lambda = \lambda_0, \mu = \mu_0$ of the matrix
of endomorphisms of $\mathcal{L}$ that defines
the eliminant has the column vector
$(v_{(\lambda_0, \mu_0)}, \ldots, D_q^{m+n-1}
v_{(\lambda_0, \mu_0)})$ in its kernel. Hence
$\Delta_{P,Q} (M, \lambda_0, \mu_0) \,
v_{(\lambda_0, \mu_0)} = 0$. Now it does not
follow automatically from this that $\Delta_{P,Q}
(M, \lambda_0, \mu_0) =0$ in $\mathcal{H}_q$
since a polynomial in $M$ might have non-trivial
kernel, as the example $(M-1) \sum_n t^n =0$
shows. However, embedding $K$ in an algebraically
closed field if necessary, we are able to show
that there exist infinitely many such
$(\lambda_0, \mu_0)$ where we \emph{can} conclude
that $\Delta_{P,Q} (M, \lambda_0, \mu_0) =0$ in
$\mathcal{H}_q$. Thus the operators $\delta_i (P,
Q)$ have an infinite-dimensional kernel, and it
is possible to show that this implies that
$\delta_i (P, Q) =0$ in $\mathcal{H}_q$.

It is more complicated to show that there exist
infinitely many $(\lambda_0, \mu_0)$ with the
required property.  The idea is to exploit the
fact that $v_{(\lambda_0, \mu_0)}$ is both in the
kernel of the operator $P - \lambda_0$ of order
$m \geq 1$ and of the operator $\Delta_{P,Q} ( M,
\lambda_0, \mu_0)$ which, if it is not zero, is
not constant. We can describe the kernel of a
non-constant polynomial element $p (M)$ of
$\mathcal{H}_q$ and the action of $P- \lambda_0$
on it explicitly enough to show that any such
$v_{(\lambda_0, \mu_0)}$ is in a subspace of
\emph{finite} dimension which (for fixed $q$)
depends only on the leading coefficient of $P$
and the degree of $P$, but not on $\lambda_0$ or
$\mu_0$. Hence for the infinity of pairs
$(\lambda_0, \mu_0)$ that can be shown to exist,
it can, by linear independence, only for finitely
many pairs be the case that $\Delta_{P,Q} (M,
\lambda_0, \mu_0)$ is not zero. For the remaining
pairs, the specialized eliminant must be zero.

\subsection{Hellstr{\"o}m-Silvestrov proof of
algebraic dependence in $\mathcal{H}_q$}

In this section we review some ideas of the proof
of Theorems \ref{thm:qBCHS} and
\ref{thm:qBCHSqfreetype} obtained by
Hellstr{\"o}m and Silvestrov in \cite{HSbook}.

Let us first $q$ be of free type. Then
$\mathrm{Z}(\mathcal{H}_q)$ is isomorphic to $K$,
and hence our aim is to prove Theorem
\ref{thm:qBCHSqfreetype}. There are three cases
to consider. In the simplest case, when $\alpha$
is of the form $cI$ for some \(c\in K\), the
polynomial \(P(x,y)=x-c\) satisfies is
annihilating for $(\alpha, \beta)$ since
\(\alpha-cI=0\).

In the second case assume that \(\alpha,\beta \)
are linear combinations of monomials with equal
degrees in $A$ and $B$ (denoted by \(\alpha,\beta
\sqsubseteq K_0\) as in \cite{HSbook} ), and that
  there is no \(c\in K\) such that \(\alpha=cI\).
  Let \(a=\deg\alpha>0\)
  and \(b=\deg\beta\). A general expression for \(P(\alpha,\beta)\), where
  $P$ has at most degree $b$ in the first variable and at most degree $a$
  in the second is
  \begin{equation} \label{Eq1: Sats 2 om algebraiskt beroende}
    \sum_{\substack{0\leqslant i\leqslant b\\0\leqslant j\leqslant a}}
p_{ij} \alpha^i \beta^j
    \text{.}
  \end{equation}
  This sum is a linear combination of the \((a+1)(b+1)\) vectors
  \( \{\alpha^i\beta^j\}_{i=0;j=0}^{b;a} \), and all these vectors
  belong to the vector space
  \(\Cen\bigl( 2ab, 0, \alpha \bigr)=
  \{\beta \sqsubseteq K_0 \mid [\alpha, \beta]=0
  \text{ and } \deg \beta \leq 2ab \}\). Now the point is that it can be shown that the dimension of this space is strictly smaller
  than \((a+1)(b+1)\). Hence there exist numbers
  \(\{p_{ij}\}_{i,j}\), not all
  zero, which make \eqref{Eq1: Sats 2 om algebraiskt beroende} zero.
  Thus there exists a $P$ as required.

  The similar algorithmic dimension growth
  type proof in the
  case when \(\alpha \not\sqsubseteq K_0\) or
  \(\beta \not\sqsubseteq K_0\) is the most
  technical part. It is based on application of
  \cite[Theorem 7.3]{HSbook} and
  \cite[Corollary 6.9]{HSbook} which are proved
  using other concepts and results in the book.
  Thus we refer the reader to \cite{HSbook} for
  details.

Let us now assume that $q$ is of torsion type and
that $p$ is its order. Observe that
$\mathcal{H}_q$, seen as a
$\mathrm{Z}(\mathcal{H}_q)$-module, contains a
spanning set of $p^2$ elements, namely
\(\{B^iA^j\}_{0 \leqslant i,j < p}\). Also
observe that $\mathrm{Z}(\mathcal{H}_q)$, by
\cite[Theorem 4.9]{HSbook} and \cite[Theorem
5.7]{HSbook}, is an integral domain. Hence any
subset of $\mathcal{H}_q$ that contains more
elements than a spanning set will be linearly
dependent. Consider the polynomial
  $$
    P(x,y) = \sum_{i=0}^p \sum_{j=0}^p c_{ij} x^i y^j
  $$
where $c_{ij} \in \mathrm{Z}(\mathcal{H}_q)$.
Clearly \(P(\alpha,\beta)\) will be a linear
combination of \((p+1)^2 > p^2\) elements in
$\mathcal{H}_q$ and these elements are linearly
dependent. Thus there are coefficients
\(\{c_{ij}\}_{0 \leqslant i,j \leqslant p}
\subset \mathrm{Z}(\mathcal{H}_q)\), not all
zero, such that \(P(\alpha,\beta)=0\). This is
exactly what the theorem claims.


\begin{thebibliography}{99}
\bibitem{Amitsur}
  S.A.\ Amitsur,
  \emph{Commutative linear differential operators},
  Pacific J.\ Math. {\bf 8} (1958), 1--10.

\bibitem{Burchnall&Chaundy1} J.L.\ Burchnall, T.W.\ Chaundy,
   \emph{Commutative ordinary differential operators},
   Proc.\ London Math.\ Soc.\ (Ser. 2) \textbf{21} (1922), 420--440.

\bibitem{Burchnall&Chaundy2}
   \bysame
   \emph{Commutative ordinary differential operators},
   Proc.\ Roy.\ Soc.\ London A \textbf{118} (1928), 557--583.

\bibitem{Burchnall&Chaundy3}
   \bysame
   \emph{Commutative ordinary differential operators. II. ---
   The Identity $P^n=Q^m$},
   Proc.\ Roy.\ Soc.\ London A \textbf{134} (1932),
   471--485.

\bibitem{HDiscMath}
L.\ Hellstr\"om, \emph{Algebraic dependence of commuting differential
operators}, Disc.\ Math.\ {\bf 231} (2001), no.\ 1--3, 246--252.

\bibitem{HSbook} L.\ Hellstr{\"o}m, S.D.\ Silvestrov,
Commuting elements in $q$-deformed Heisenberg algebras, World Scientific, New
Jersey, 2000.

\bibitem{LHSSGWSergpotJAlg}
L.\ Hellstr{\"o}m, S.\ Silvestrov, \emph{Ergodipotent maps and commutativity of
elements in non-commutative rings and algebras with twisted intertwining}, J.\
Algebra {\bf 314} (2007), 17-41.

\bibitem{JSSvBCqHeis} M.\ de Jeu, P.C.\ Svensson, S.\ Silvestrov, \emph{Algebraic curves for commuting
elements in the $q$-deformed Heisenberg algebra}, arXiv:0710.2748v1 [math.RA],
2007, 17pp., to appear.

\bibitem{KrichIntro1} I.M.\ Krichever, \emph{Integration of non-linear equations by the
methods of algebraic geometry}, Funktz.\ Anal.\ Priloz.\ \textbf{11}, 1 (1977),
15--31.

\bibitem{KrichIntro2} I.M.\ Krichever, \emph{Methods of algebraic geometry in the theory
of nonlinear equations}, Uspekhi\ Mat.\ Nauk,
\textbf{32}, 6 (1977), 183--208.

\bibitem{LarssonSilv} D.\ Larsson, S.D.\ Silvestrov,
\emph{Burchnall-Chaundy theory for $q$-difference operators and $q$-deformed
Heisenberg algebras}, J.\ Nonlin.\ Math.\ Phys.\ {\bf 10} (2003), suppl.\ 2,
95-106.

\bibitem{MumIntro} D.\ Mumford, \emph{An algebro-geometric construction of commuting operators and
   of solutions to the Toda lattice equation, Korteweg--de~Vries
   equation and related non-linear equations},
   Proc.\ Int.\ Symp.\ on Algebraic Geometry, Kyoto (1978), 115--153.

\bibitem{Nesterenko} Yu.\ Nesterenko,
{\it Modular functions and transcendence problems}, C.R.\ Acad.\ Sci.\ Paris
Ser.\ I Math.\ {\bf 322} (1996), no.\ 10, 909--914.

\bibitem{NestPhil} Yu.\ Nesterenko, P.\ Philippon,
(Eds.), Introduction to algebraic independence theory, Lecture Notes in
Mathematics 1752. Springer-Verlag, Berlin, 2001.

\end{thebibliography}
\end{document}